\documentstyle[leqno,11pt]{article}

\newtheorem{th}{Theorem}[section]
\newtheorem{prop}[th]{Proposition}
\newtheorem{cor}[th]{Corollary}
\newcounter{defin}[section]
\renewcommand{\thedefin}{\thesection.\arabic{defin}}
\newcounter{ex}[section]
\renewcommand{\theex}{\thesection.\arabic{ex}}
\newcounter{rem}[section]
\renewcommand{\therem}{\thesection.\arabic{rem}}
\title{Generalized Metrical Multi-Time Lagrange\\
Geometry of Physical Fields}
\author{Mircea Neagu}
\date{}
\begin{document}

\maketitle
\begin{abstract}
Section 1 contains some physical and geometrical aspects of the generalized
Lagrangian geometry of physical fields developed by Miron and Anastasiei
\cite{7}, which represents the start point in our generalized metrical
multi-time Lagrangian approach of the theory of physical fields. Section 2
exposes a geometrization of a Kronecker $h$-regular vertical  metrical d-tensor
$G^{(\alpha)(\beta)}_{(i)(j)}$ on the jet fibre bundle of order one $J^1(T,M)$.
We emphasize that this geometrization gives a mathematical model for both the
gravitational and electromagnetic field theory, in a general setting. Thus,
Section 3 presents the generalized metrical multi-time Lagrange theory of
electromagnetism and describes its Maxwell equations. Section 4 presents the
Einstein equations which govern the generalized  metrical multi-time Lagrange
theory of gravitational field. The conservation laws of the gravitational field
are also described in terms of generalized metrical multi-time Lagrange geometry.
Section 5 gives a more natural form of the Einstein equations and of their
conservation laws.
\end{abstract}

\section{Generalized Lagrangian theory
of physical fields}

\hspace{5mm}  About ten years ago, Miron and Anastasiei were  developed the
{\it generalized Lagrangian theory of physical fields} \cite{7}, which naturally
generalizes the Finslerian and Lagrangian ones. The geometrical background
of this theory relies on the notion of {\it generalized Lagrange space}
$GL^n=(M,g_{ij}(x^k,y^k))$, which consists of a real, $n$-dimensional manifold
$M$ coordinated by $(x^i)_{i=\overline{1,n}}$, and a symmetric, of  rank  $n$
and of constant signature {\it fundamental metrical d-tensor $g_{ij}(x^k,y^k)$}
on $TM$.

In the sequel, we try to expose the main geometrical and physical aspects of
this  field theory.

From physical point of view, the fundamental metrical d-tensor has the physical meaning of an
{\it "unified"} gravitational field on $TM$, which consists of one {\it
"external"} $(x)$-gravitational field spanned by points $\{x\}$, and the other
{\it "internal"} $(y)$-gravitational field spanned by directions $\{y\}$. It
should be emphasized that $y$ is endowed with some microscopic character of
the space-time structure. Moreover, since $y$ is a vector field different of
an ordinary vector field, the $y$-dependence has combined with the concept of
{\it anisotropy}.

The field theory developed on a generalized Lagrange space $GL^n$ relies on
a fixed {\it "a priori"} nonlinear connection $\Gamma=(N^i_j(x,y))$  on the
tangent bundle $TM$. This plays the role of mapping operator of internal
$(y)$-field on the external $(x)$-field, and prescribes the {\it "interaction"} between $(x)$-
and $(y)$- fields. From geometrical point of view, the nonlinear connection
allows the construction of the {\it adapted bases}
\begin{equation}
\begin{array}{l}\medskip
\displaystyle{\left\{{\delta\over\delta x^i}={\partial\over\partial x^i}-
N_i^j{\partial\over\partial y^j},{\partial\over\partial y^i}\right\}\subset
{\cal X}(TM)},\\
\{dx^i,\delta y^i=dy^i+N^i_jdx^j\}\subset{\cal X}^*(TM).
\end{array}
\end{equation}

As to the spatial structure, the most important thing is to determine the {\it
Cartan canonical connection} $C\Gamma=(L^i_{jk},C^i_{jk})$ with respect to $g_{ij}$,
which  comes from the metrical  conditions
\begin{equation}
\left\{\begin{array}{l}\medskip
\displaystyle{g_{ij\vert k}={\delta g_{ij}\over\delta x^k}-L^m_{ik}g_{mj}-
L^m_{jk}g_{mi}=0}\\
\displaystyle{g_{ij}\vert_k={\partial g_{ij}\over\partial y^k}-C^m_{ik}g_{mj}-
C^m_{jk}g_{mi}=0},
\end{array}\right.
\end{equation}
where $"_{\vert k}"$ and $"\vert_k"$  are  the local $h$-  and $v$- covariant
derivatives of  $C\Gamma$. The importance  to the Cartan canonical connection
comes from its main role played in the generalized Lagrangian theory of physical
fields.

Concerning the {\it "unified"} field $g_{ij}(x,y)$ of $GL^n$, the authors
constructed a Sasakian-like metric on $TM$,
\begin{equation}
G=g_{ij}dx^i\otimes dx^j+g_{ij}\delta y^i\otimes\delta y^j.
\end{equation}

In this context, the Einstein equations of the gravitational potentials
$g_{ij}(x,y)$ of a generalized Lagrange space $GL^n,\;n>2$, are postulated as being the
Einstein equations attached to $C\Gamma$ and $G$, namely \cite{7},
\begin{equation}
\left\{\begin{array}{ll}\medskip
\displaystyle{R_{ij}-{1\over 2}Rg_{ij}={\cal K}{\cal T}_{ij}^H},&^\prime P_{ij}
={\cal  K}{\cal T}_{ij}^1,\\
\displaystyle{S_{ij}-{1\over 2}Sg_{ij}={\cal K}{\cal T}_{ij}^V},&^{\prime\prime}
P_{ij}=-{\cal  K}{\cal T}_{ij}^2,
\end{array}\right.
\end{equation}
where $R_{ij}=R^m_{ijm}$, $S_{ij}=S^m_{ijm}$, $^\prime P_{ij}=P^m_{ijm}$,
$^{\prime\prime} P_{ij}=P^m_{imj}$ are the Ricci tensors of $C\Gamma$, $R=
g^{ij}R_{ij}$, $S=g^{ij}S_{ij}$ are the scalar curvatures, ${\cal T}_{ij}^H$,
${\cal T}_{ij}^V$, ${\cal T}_{ij}^1$, ${\cal T}_{ij}^2$ are the components of the
energy-momentum tensor ${\cal T}$ and ${\cal K}$ is the Einstein constant (equal to 0 for
vacuum). Moreover, the energy-momentum tensors ${\cal T}_{ij}^H$ and ${\cal T}
_{ij}^V$ satisfy the so-called {\it conservation laws}
\begin{equation}
{\cal K}{\cal T}^{H\;m}_{\quad j\vert m}=-{1\over 2}(P^{hm}_{js}R^s_{hm}+2R^s_{mj}
P^m_s),\quad{\cal K}{\cal T}^{V\;m}_{\quad j\vert m}=0,
\end{equation}
where all notations are described in \cite{7}.

The generalized Lagrangian theory of electromagnetism relies on the {\it canonical Liouville
vector field} {\bf C}=$\displaystyle{y^i{\partial\over\partial y^i}}$ and the Cartan canonical
connection $C\Gamma$ of the generalized Lagrange space $GL^n$. In this context, the authors
introduce the {\it electromagnetic 2-form} on $TM$,
\begin{equation}
F=F_{ij}\delta y^i\wedge  dx^j+f_{ij}\delta y^i\wedge\delta y^j,
\end{equation}
where
\begin{equation}
\begin{array}{l}\medskip
\displaystyle{F_{ij}={1\over 2}[(g_{im}y^m)_{\vert j}-(g_{jm}y^m)_{\vert i}],}\\
\displaystyle{f_{ij}={1\over 2}[(g_{im}y^m)\vert_j-(g_{jm}y^m)\vert_i].}
\end{array}
\end{equation}

Using the Bianchi identities attached to the Cartan canonical
connection $C\Gamma$, they conclude that the electromagnetic  components
$F_{ij}$ and $f_{ij}$ are governed by the following {\it equations of Maxwell
type},
\begin{equation}
\hspace*{5mm}\left\{\begin{array}{l}\medskip
F_{ij\vert k}+F_{jk\vert i}+F_{ki\vert j}=-\sum_{\{i,j,k\}}\left[C_{imr}y^m+
(g_{im}y^m)_{\vert r}\right]R^r_{jk}\\\medskip
F_{ij}\vert_k+F_{jk}\vert_i+F_{ki}\vert_j=-(f_{ij\vert k}+f_{jk\vert i}+
f_{ki\vert j})\\
f_{ij}\vert_k+f_{jk}\vert_i+f_{ki}\vert_j=0.
\end{array}\right.
\end{equation}
For a deeply exposition of the classical generalized Lagrange theory of
physical fields, the reader is invited to see \cite{7}.

In this paper, we naturally extend the
previous field theory to a general one, constructed on the  jet fibre bundle
of order one $J^1(T,M)\to T\times M$, where $T$ is  a smooth, real,  $p$-dimensional
{\it "multi-time"} manifold  coordinated  by $(t^\alpha)_{\alpha=\overline{1,p}}$
and $M$ is a smooth, real $n$-dimensional {\it "spatial"} manifold coordinated
by $(x^i)_{i=\overline{1,n}}$. The gauge group of $J^1(T,M)$ is
\begin{equation}\label{gg}
\left\{\begin{array}{l}\medskip
\tilde t^\alpha=\tilde t^\alpha(t^\beta)\\\medskip
\tilde x^i=\tilde x^i(x^j)\\
\displaystyle{\tilde x^i_\alpha={\partial\tilde x^i\over\partial x^j}
{\partial t^\beta\over\partial\tilde t^\alpha}x^j_\beta},
\end{array}\right.
\end{equation}
where the meaning of $x^i_\alpha$  is that of {\it partial derivatives} or,
alternatively, of {\it  partial directions}.
Consequently, our field theory can be regarded
as a theory in which  the physical entities  are  dependent of some temporal
coordinates  $(t^\alpha)$, spatial coordinates  $(x^i)$ and partial directions
$(x^i_\alpha)$.
Our field theory is created, in a natural manner, from a given {\it Kronecker
$h$-regular vertical fundamental metrical d-tensor} $G^{(\alpha)(\beta)}_{(i)
(j)}$ on $J^1(T,M)$, and can be called the {\it generalized metrical multi-time Lagrange
theory of physical fields.}

Finally, we recall that the jet fibre bundle of order one  is a basic object in
the study of classical and quantum field theories.

\section{Generalized metrical multi-time Lagrange spaces}

\setcounter{equation}{0}
\hspace{5mm}  Let us consider $T$ (resp. $M$) a {\it "temporal"}  (resp. {\it
"spatial"}) manifold of dimension $p$ (resp. $n$), coordinated by $(t^\alpha)_
{\alpha=\overline{1,p}}$ (resp. $(x^i)_{i=\overline{1,n}}$). Let $E=J^1(T,M)\to
T\times M$ be the jet fibre bundle of order one associated to these manifolds.
The {\it bundle of configurations} $J^1(T,M)$ is coordinated by
$(t^\alpha,x^i,x^i_\alpha)$, where $\alpha=\overline{1,p}$ and $i=\overline{1,n}$.
Note that, throughout  this paper, the indices $\alpha,\beta,\gamma,
\ldots$ run from $1$ to $p$, and the indices $i,j,k,\ldots$ run from $1$ to $n$.
\medskip\\
\addtocounter{defin}{1}
{\bf Definition \thedefin} A d-tensor field on $E$, defined by the local
components $G^{(\alpha)(\beta)}_{(i)(j)}(t^\gamma,x^k,x^k_\gamma)$, is called
a {\it vertical fundamental metrical d-tensor}.\medskip\\
\addtocounter{ex}{1}
{\bf Example \theex} Let $L:E\to R$ be a  multi-time Lagrangian function. Then,
the local components
\begin{equation}\label{vfm}
G^{(\alpha)(\beta)}_{(i)(j)}(t^\gamma,x^k,x^k_\gamma)={1\over 2}{\partial^2 L
\over\partial x^i_\alpha\partial x^j_\beta}
\end{equation}
represent a vertical fundamental metrical d-tensor on  $E$, which is called the
{\it canonical vertical fundamental metrical d-tensor attached to the Lagrangian
function $L$}.\medskip

Now, let $h=(h_{\alpha\beta})$ be a fixed semi-Riemannian metric on the temporal
manifold $T$ and $g_{ij}(t^\gamma, x^k, x^k_\gamma)$ be a d-tensor on $E$, symmetric,
of rank $n$, and having a constant signature.\medskip\\
\addtocounter{defin}{1}
{\bf Definition \thedefin} A vertical fundamental metrical d-tensor
$G^{(\alpha)(\beta)}_{(i)(j)}$ which is
of the form
\begin{equation}
G^{(\alpha)(\beta)}_{(i)(j)}(t^\gamma,x^k,x^k_\gamma)=h^{\alpha\beta}(t^
\gamma)g_{ij}(t^\gamma,x^k,x^k_\gamma),
\end{equation}
is called a {\it Kronecker $h$-regular vertical fundamental metrical d-tensor
with respect to the temporal semi-Riemannian metric $h=(h_{\alpha\beta})$}.
\medskip\\
\addtocounter{rem}{1}
{\bf Remark \therem} In the case $\dim T\geq 2$, the condition of Kronecker
$h$-regularity imposed to the canonical vertical fundamental metrical d-tensor of a
multi-time Lagrangian function $L$, i. e.
\begin{equation}
\displaystyle{G^{(\alpha)(\beta)}_{(i)(j)}(t^\gamma,
x^k,x^k_\gamma)={1\over 2}{\partial^2 L\over\partial x^i_\alpha\partial x^j_
\beta}=h^{\alpha\beta}(t^\gamma)g_{ij}(t^\gamma,x^k,x^k_\gamma)},
\end{equation}
implies the independence of partial directions of the d-tensor field $g_{ij}$,
that is, $g_{ij}(t^\gamma,x^k,x^k_\gamma)=g_{ij}(t^\gamma,x^k)$. In other words,
a {\it Kronecker $h$-regular multi-time Lagrange function}  take the form
\begin{equation}
L=h^{\alpha\beta}(t)g_{ij}(t,x)x^i_\alpha x^j_\beta+U^{(\alpha)}_{(i)}(t,x)x^i_
\alpha+F(t,x).
\end{equation}
For more details, see the characterization theorem of metrical multi-time
Lagrange spaces, described in \cite{13}.
\medskip\\
\addtocounter{ex}{1}
{\bf Example \theex} Let $G^{(\alpha)(\beta)}_{(i)(j)}(t^\gamma,x^k,x^k_
\gamma)$ be the canonical vertical fundamental metrical d-tensor of the multi-time
Lagrangian function  $L$, not necessarily a Kronecker $h$-regular one. Let
us suppose that the contraction d-tensor
\begin{equation}
g_{ij}(t^\gamma,x^k,x^k_\gamma)={1\over p}h_{\mu\nu}(t^\gamma)
G^{(\mu)(\nu)}_{(i)(j)}(t^\gamma,x^k,x^k_\gamma),
\end{equation}
where $p=\dim T$, is symmetric, of rank $n$ and having a constant signature
on $E$. Then, the vertical fundamental metrical d-tensor
\begin{equation}
{\cal G}^{(\alpha)(\beta)}_{(i)(j)}(t^\gamma,x^k,x^k_\gamma)=h^{\alpha\beta}(t^
\gamma)g_{ij}(t^\gamma,x^k,x^k_\gamma),
\end{equation}
is a Kronecker  $h$-regular one. This  is called the {\it canonical Kronecker
$h$-regular vertical fundamental metrical d-tensor associated to the
Lagrangian function $L$}.\medskip\\
\addtocounter{rem}{1}
{\bf Remark \therem} If $L$ is a Kronecker $h$-regular multi-time Lagrangian
function, the equality ${\cal G}^{(\alpha)(\beta)}_{(i)(j)}=G^{(\alpha)
(\beta)}_{(i)(j)}$ holds good.\medskip

In this context, we can introduce the  following\medskip\\
\addtocounter{defin}{1}
{\bf  Definition \thedefin} If $p=\dim T$ and $n=\dim M$, a pair
$$
GML^n_p=(J^1(T,M),G^{(\alpha)(\beta)}_{(i)(j)}),
$$
which consists of the
1-jet fibre bundle and a Kronecker $h$-regular vertical fundamental metrical
d-tensor $G^{(\alpha)(\beta)}_{(i)(j)}(t^\gamma,x^k,x^k_\gamma)=h^{\alpha\beta}
(t^\gamma)g_{ij}(t^\gamma,x^k,x^k_\gamma)$, is called a {\it  generalized
metrical multi-time Lagrange space}.\medskip\\
\addtocounter{rem}{1}
{\bf Remarks \therem} i) Let us consider the particular case $T=R$ of the usual time
axis  represented by the set of real numbers.  In this  case, the  coordinates
of $J^1(R,M)\equiv R\times TM$ are denoted by $(t,x^i,y^i)$.
From physical point of view, we emphasize that the fibre bundle
\begin{equation}
J^1(R,M)\to R\times M,\quad (t,x^i,y^i)\to (t,x^i),
\end{equation}
is regarded like a {\it bundle of configurations}, in mechanics terms.
The gauge group of this bundle of configurations is
\begin{equation}\label{G_1}
\left\{\begin{array}{l}
\tilde t=\tilde t(t)\\
\tilde x^i=\tilde x^i(x^j)\\
\displaystyle{\tilde y^i={\partial\tilde x^i\over\partial x^j}{dt\over
d\tilde t}y^j.}
\end{array}\right.
\end{equation}

We  remark that the form of this gauge group stands out by the {\it relativistic}
character of the time $t$. For that reason, in the particular case $(T,h)=(R,\delta)$,
a generalized metrical multi-time Lagrange space is called a {\it  generalized
relativistic rheonomic Lagrange space} and is denoted $GRL^n=(J^1(R,M),
G^{(1)(1)}_{(i)(j)}(t,x^k,y^k)=g_{ij}(t,x^k,y^k))$.

It is important to note that, in the {\it classical generalized rheonomic
Lagrangian geometry} \cite{7}, the bundle of configuration is the fibre
bundle
\begin{equation}
R\times TM\to M,\;(t,x^i,y^i)\to (x^i),
\end{equation}
whose geometrical invariance group is
\begin{equation}\label{G_2}
\left\{\begin{array}{l}
\tilde t=t\\
\tilde x^i=\tilde x^i(x^j)\\
\displaystyle{\tilde y^i={\partial\tilde x^i\over\partial x^j}y^j.}
\end{array}\right.
\end{equation}
Obviously, the structure of the gauge group \ref{G_2} emphasizes the  {\it
absolute}  character of the time $t$ from the classical generalized rheonomic
Lagrangian geometry. At the same time, we point out that the gauge group \ref{G_2}
is a subgroup of \ref{G_1}. In other words, the gauge group of the jet bundle of
order one from the generalized relativistic rheonomic Lagrangian geometry is
more general than that used in the classical generalized rheonomic Lagrangian
geometry, which ignores the temporal reparametrizations.

Finally, we invite the  reader to  compare our generalized rheonomic Lagrange
geometry developed on  $J^1(R,M)$ to that sketched by Miron and Anastasiei on
$R\times TM$, at the end of the book \cite{7}.

ii) If  the temporal manifold $T$ is 1-dimensional, then,
via a temporal reparametrization, we have  $J^1(T,M)\equiv J^1(R,M)$.
In other words, a generalized  metrical multi-time Lagrangian space having
$\dim T=1$ is a  {\it reparametrized generalized relativistic rheonomic
Lagrange space}.\medskip\\
\addtocounter{ex}{1}
{\bf Example \theex} Let $U^{(\alpha)}_{(i)}(t^\gamma,x^k)$ be a d-tensor on
the jet space  $J^1(T,M)$
and $F:T\times M\to R$ be a smooth function.  Let us consider\linebreak
$L:J^1(T,M)\to
R,$ the quadratic
multi-time Lagrangian function used in  various physical models,
\begin{equation}
\hspace*{5mm}L=G^{(\alpha)(\beta)}_{(i)(j)}(t^\gamma,x^k)x^i_\alpha
x^j_\beta+U^{(\alpha)}_{(i)}(t^\gamma,x^k)x^i_\alpha+F(t^\gamma,x^k),
\end{equation}
whose vertical fundamental metrical d-tensor $G^{(\alpha)(\beta)}_{(i)(j)}
(t^\gamma,x^k)$ is symmetric, of rank $n$ and has a constant signature with respect to
the indices $(i)$ and $(j)$. In these conditions, taking {\it the canonical
Kronecker $h$-regular vertical metrical d-tensor of $L$},
\begin{equation}
{\cal G}^{(\alpha)(\beta)}_{(i)(j)}(t^\gamma,x^k)={1\over p}h^{\alpha
\beta}(t^\gamma)h_{\mu\nu}(t^\gamma)G^{(\mu)(\nu)}_{(i)(j)}(t^\gamma,x^k),
\end{equation}
where $p=\dim T$,  the pair
$
GML^n_p=(J^1(T,M),{\cal G}^{(\alpha)(\beta)}_{(i)(j)}(t^\gamma,x^k))
$
is a  generalized metrical multi-time Lagrange space. This is called the {\it
canonical generalized multi-time Lagrange space attached to the quadratic multi-time
Lagrangian function  $L$}.\medskip\\
\addtocounter{ex}{1}
{\bf Example \theex} More  general, let us consider an arbitrary Lagrangian
function
$L:J^1(T,M)\to R$, whose vertical fundamental metrical d-tensor
$$
G^{(\alpha)(\beta)}_{(i)(j)}(t^\gamma,x^k,x^k_\gamma)={1\over 2}{\partial^2 L
\over\partial x^i_\alpha\partial x^j_\beta}
$$
is symmetric, of rank $n$ and has a constant signature with respect to the
indices $(i)$ and $(j)$. Then, the  pair
\begin{equation}
GML^n_p=(J^1(T,M),{\cal G}^{(\alpha)(\beta)}_{(i)(j)}(t^\gamma,x^k,x^k_\gamma)),
\end{equation}
where  ${\cal G}^{(\alpha)(\beta)}_{(i)(j)}(t^\gamma,x^k,x^k_\gamma)$ is the
canonical Kronecker $h$-regular vertical metrical d-tensor of $L$, is called
the {\it canonical generalized multi-time Lagrange space attached to the
multi-time Lagrangian function  $L$}.\medskip\\
\addtocounter{ex}{1}
{\bf Example \theex} Let $\varphi_{ij}(x^k)$  be  a semi-Riemannian metric on
the spatial manifold  $M$ and let $\sigma:J^1(T,M)\to R$ be a conformal smooth
function,  which gives the {\it magnitude} of directions $x^i_\alpha$. In this context,
the pair $GML^n_p=(J^1(T,M),G^{(\alpha)(\beta)}_{(i)(j)}(t^\gamma,x^k,x^k_
\gamma))$, where
\begin{equation}
G^{(\alpha)(\beta)}_{(i)(j)}(t^\gamma,x^k,x^k_\gamma))=h^{\alpha\beta}(t^\gamma)
e^{2\sigma(t^\gamma,x^k, x^k_\gamma)}\varphi_{ij}(x^k),
\end{equation}
is a  generalized  metrical multi-time Lagrange space. Note that, in this case,
we have
\begin{equation}
g_{ij}(t^\gamma,x^k,x^k_\gamma)=e^{2\sigma(t^\gamma,x^k, x^k_\gamma)}
\varphi_{ij}(x^k).
\end{equation}
From  physical point of view,  the  interesting properties
of this space are obtained  considering the special conformal  functions:
\medskip

i) $\sigma=U^{(\alpha)}_{(i)}(t^\gamma,x^k)x^i_\alpha$,\medskip

ii) $\sigma=h^{\alpha\beta}(t^\gamma)A_i(x^k)A_j(x^k)x^i_\alpha x^j_\beta$,
\medskip

iii) $\sigma=\varphi_{ij}(x^k)X^\alpha(t^\gamma)X^\beta(t^\gamma)x^i_\alpha
x^j_\beta$,\medskip\\
where $U^{(\alpha)}_{(i)}(t^\gamma,x^k)$ is a d-tensor on $E$, $A_i(x^k)$ is
a covector field on $M$, and $X^\alpha(t^\gamma)$ is a vector field on $T$. A
deeply geometrical and  physical study of this generalized metrical multi-time
Lagrange space is made in \cite{10}.
\medskip\\
\addtocounter{ex}{1}
{\bf Example \theex}  Let us consider the Kronecker  $h$-regular vertical
metrical d-tensor
\begin{equation}
G^{(\alpha)(\beta)}_{(i)(j)}=h^{\alpha\beta}(t^\gamma)
\left[\varphi_{ij}(x^k)+\left(1-{1\over n(t^\gamma,x^k,x^k_\gamma)}\right)
Y_iY_j\right],
\end{equation}
where $n:J^1(T,M)\to [1,\infty)$ is a  smooth function  representing the
{\it refraction indices of medium}, $X^\alpha(t^\gamma)$ are  the components
of a  vector field on $T$ representing  the {\it direction of refraction},
and
$$
Y_i(t^\gamma,x^k,x^k_\gamma)=\varphi_{im}(x^k)x^m_\mu X^\mu(t^\gamma).
$$
The generalized metrical multi-time Lagrange space
$GML^n_p=(J^1(T,M),G^{(\alpha)(\beta)}_{(i)(j)}(t^\gamma,x^k,x^k_\gamma))$ is
called  the  {\it generalized metrical multi-time Lagrange  space  of relativistic
geometrical optic}. It is obvious that, in this case, we have
\begin{equation}
g_{ij}(t^\gamma,x^k,x^k_\gamma)=\varphi_{ij}(x^k)+\left(1-{1\over n(t^\gamma,
x^k,x^k_\gamma)}\right)Y_iY_j.
\end{equation}
Moreover, by a direct calculation, it follows
\begin{equation}\hspace*{5mm}
g^{ij}(t^\gamma,x^k,x^k_\gamma)=\varphi^{ij}(x^k)+{1-\displaystyle{{1\over n
(t^\gamma,x^k,x^k_\gamma)}}\over 1+\displaystyle{\left(1-{1\over n(t^\gamma,
x^k,x^k_\gamma)}\right)Y^2}}Y^iY^j,
\end{equation}
where $Y^i=\varphi^{ir}Y_r$ and $Y^2=Y^mY_m$.\medskip\\
\addtocounter{rem}{1}
{\bf Remarks \therem} i) The terminology used above is a natural one, because
the metrical multi-time Lagrange  space of relativistic geometrical optic is
a  natural generalization of that so-called the {\it generalized Lagrange space
of relativistic geometrical  optic} from Miron-Anastasiei theory \cite{7}.

ii) The geometry  of  the generalized metrical multi-time Lagrange  space of
relativistic geometrical optic will be deeply  studied in \cite{9}.\medskip

In order to develope a geometry on a generalized metrical multi-time Lagrange
space  $GML^n_p$, whose Kronecker $h$-regular vertical fundamental metrical
d-tensor is
\begin{equation}
G^{(\alpha)(\beta)}_{(i)(j)}(t^\gamma,x^k,x^k_\gamma)=h^{\alpha
\beta}(t^\gamma)g_{ij}(t^\gamma,x^k,x^k_\gamma),
\end{equation}
we need  a  nonlinear connection $\Gamma=(M^{(i)}_{(\alpha)\beta},N^{(i)}_{(
\alpha)j})$ on $E=J^1(T,M)$.

Firstly, we point out that the vertical fundamental metrical d-tensor
$G^{(\alpha)(\beta)}_{(i)(j)}(t^\gamma,x^k,x^k_\gamma)$ of $GML^n_p$ induces
a {\it canonical temporal nonlinear connection} \cite{12}, defined by the local
components
\begin{equation}
M^{(i)}_{(\alpha)\beta}=-H^\gamma_{\alpha\beta}x^i_\gamma,
\end{equation}
where  $H^\gamma_{\alpha\beta}$ are  the Christoffel symbols of the temporal
semi-Riemannian metric $h_{\alpha\beta}$.

Secondly, to associate  a {\it canonical spatial nonlinear connection} $N^{(i)}_
{(\alpha)j}$ to the vertical fundamental metrical d-tensor $G^{(\alpha)(\beta)}
_{(i)(j)}(t^\gamma,x^k,x^k_\gamma)$ of $GML^n_p$, we introduce the concept of
{\it absolut energy  Lagrangian function ${\cal  E}:J^1(T,M)\to R,$ attached
to a $GML^n_p$}, namely,
\begin{equation}
{\cal E}=G^{(\mu)(\nu)}_{(m)(r)}x^m_\mu x^r_\nu=h^{\mu\nu}(t^\gamma)g_{mr}(t^\gamma,x^k,x^k_
\gamma)x^m_\mu x^r_\nu.
\end{equation}

{\bf 1.} Let us suppose that ${\cal E}$ is a Kronecker $h$-regular Lagrangian function
\cite{13}. In this context,  following the development of the metrical multi-time
Lagrange geometry \cite{13}, we can construct a canonical spatial nonlinear
connection $N^{(i)}_{(\alpha)j}$, directly from  ${\cal E}$.  For example, let
us consider
$
GML^n_p=(J^1(T,M),{\cal G}^{(\alpha)(\beta)}_{(i)(j)}(t^\gamma,x^k)),
$
the canonical generalized multi-time Lagrange space attached to the quadratic
multi-time Lagrangian function  $L$ from the example {\bf 2.3}. In this case,
since the vertical fundamental metrical d-tensor is
\begin{equation}\label{mL}
{\cal G}^{(\alpha)(\beta)}_{(i)(j)}(t^\gamma,x^k)=h^{\alpha\beta}(t^\gamma)
g_{ij}(t^\gamma,x^k),
\end{equation}
where
$\displaystyle{g_{ij}(t^\gamma,x^k)={1\over p}h_{\mu\nu}(t^\gamma)G^{(\mu)
(\nu)}_{(i)(j)}(t^\gamma,x^k)}$, it follows that the absolute energy Lagrange
function becomes
\begin{equation}
{\cal E}={\cal G}^{(\mu)(\nu)}_{(m)(r)}x^m_\mu x^r_\nu=h^{\mu\nu}(t^\gamma)g_{mr}(t^
\gamma,x^k)x^m_\mu x^r_\nu,
\end{equation}
that is, it is a Kronecker  $h$-regular  Lagrange function. Consequently,
the canonical spatial nonlinear connection of this generalized multi-time Lagrange
space is \cite{13}
\begin{equation}\label{nlcq}
N^{(i)}_{(\alpha)j}=\Gamma^i_{jm}x^m_\alpha+{g^{im}\over 2}{\partial g_{jm}
\over\partial t^\alpha},
\end{equation}
where
\begin{equation}
\Gamma^i_{jk}(t^\mu,x^m)={g^{ir}\over 2}\left({\partial g_{jr}\over\partial
x^k}+{\partial g_{kr}\over\partial x^j}-{\partial g_{jk}\over\partial x^r}
\right)
\end{equation}
are  the {\it generalized Christoffel symbols} of the multi-time dependent
spatial metric $g_{ij}(t^\gamma,x^k)$.

{\bf 2.} If ${\cal E}$ is not a Kronecker $h$-regular Lagrangian function, we are
forced to give an {\it "a priori"}  spatial nonlinear connection. In this sense,
in the examples {\bf 2.5} and {\bf 2.6}, it is convenient to use the spatial
nonlinear connection
\begin{equation}
N^{(i)}_{(\alpha)j}=\gamma^i_{jm}x^m_\alpha,
\end{equation}
where $\gamma^i_{jm}$ are the Christoffel symbols of the spatial semi-Riemannian
metric $\varphi_{ij}$. The using of this spatial nonlinear
connection is justified, esspecialy by its physical aspects, in \cite{9}
and  \cite{10}.  Nevertheless, from a geometrical point  of view, we point
out that this spatial nonlinear connection is buildet  directly from the
vertical fundamental metrical d-tensor $G^{(\alpha)(\beta)}_{(i)(j)}$ (see
the examples {\bf 2.5} and {\bf 2.6}). Consequently, we have a {\it "metrical"}
character (see \cite{4}) of the  geometry that we will construct on these spaces.
In other words, all geometrical objects with physical meaning that we will
build in \cite{9} and \cite{10}, will be directly arised from
$G^{(\alpha)(\beta)}_{(i)(j)}$.

Finally, we consider that it is important to study what conditions must be imposed to
the absolute energy Lagrangian function ${\cal E}$, in order to obtain its
Kronecker $h$-regularity. In order to do that, let us suppose $\dim T\geq  2$.
In this context,
we proved in \cite{13} that the Kronecker $h$-regularity  of ${\cal E}$ reduces
to the existence of some d-tensors $\varphi_{ij}(t^\gamma,x^k)$, $U^{(\alpha)}_{(i)}
(t^\gamma,x^k)$ and of some  smooth function $F:T\times M\to R$, such that,
\begin{equation}
{\cal E}=h^{\alpha\beta}(t^\gamma)\varphi_{ij}(t^\gamma,x^k)x^i_\alpha
x^j_\beta+U^{(\alpha)}_{(i)}(t^\gamma,x^k)x^i_\alpha+F(t^\gamma,x^k),
\end{equation}
that is,
$$
h^{\alpha\beta}(t^\gamma)\left[g_{ij}(t^\gamma,x^k,x^k_\gamma)-\varphi_{ij}
(t^\gamma,x^k)\right]x^i_\alpha x^j_\beta-U^{(\alpha)}_{(i)}(t^\gamma,x^k)x^i_\alpha
-F(t^\gamma,x^k)=0.
$$
Consequently, if the spatial metrical d-tensor $g_{ij}$ of a $GML^n_p$ does
not depend of the partial directions  $x^i_\alpha$ (see \ref{mL}), putting
$\varphi_{ij}(t^\gamma,x^k)=g_{ij}(t^\gamma,x^k)$, $U^{(\alpha)}_{(i)}(t^
\gamma,x^k)=0$  and $F(t^\gamma,x^k)=0$, it follows that ${\cal E}$ is a
Kronecker $h$-regular Lagrangian function. In conclusion, in this case,
we can build a natural nonlinear spatial connection,  directly from the
vertical fundamental metrical d-tensor of $GML^n_p$ (see \ref{nlcq}).
\medskip\\
{\bf\underline{Open problems}. i)} Supposing $\dim T\geq 2$, are there the
spatial metrical d-tensors $g_{ij}(t^\gamma,x^k,x^k_\gamma)$, depending
effectively of $x^i_\alpha$, such that ${\cal E}$ to be a Kronecker $h$-regular
Lagrangian function?\medskip

{\bf ii)} Is it possible  to construct, in a natural way, a spatial nonlinear
connection $N^{(i)}_{(\alpha)j}$ from the vertical fundamental metrical d-tensor
$G^{(\alpha)(\beta)}_{(i)(j)}$  of a generalized metrical multi-time  Lagrange
space $GML^n_p$?\medskip

Now, let us consider a generalized metrical multi-time  Lagrange space
$$
GML^n_p=(J^1(T,M),G^{(\alpha)(\beta)}_{(i)(j)}(t^\gamma,x^k)=h^{\alpha\beta}(t^\gamma)g_{ij}
(t^\gamma,x^k,x^k_\gamma))
$$
and an {\it "a priori"} fixed nonlinear connection
$\Gamma=(M^{(i)}_{(\alpha)\beta},N^{(i)}_{(\alpha)j})$ on the jet space
$J^1(T,M)$.
Let $\displaystyle{\left\{{\delta\over\delta
t^\alpha}, {\delta\over\delta x^i}, {\partial\over\partial x^i_\alpha}\right\}
\subset{\cal X}(E)}$ and $\{dt^\alpha, dx^i, \delta x^i_\alpha\}\subset{\cal
X}^*(E)$ be the adapted bases of the nonlinear connection $\Gamma$, where
\begin{equation}
\left\{\begin{array}{l}\medskip
\displaystyle{{\delta\over\delta t^\alpha}={\partial\over\partial t^\alpha}-
M^{(j)}_{(\beta)\alpha}{\partial\over\partial x^j_\beta}}\\\medskip
\displaystyle{{\delta\over\delta x^i}={\partial\over\partial x^i}-
N^{(j)}_{(\beta)i}{\partial\over\partial x^j_\beta}}\\
\delta x^i_\alpha=dx^i_\alpha+M^{(i)}_{(\alpha)\beta}dt^\beta+N^{(i)}_{(\alpha)
j}dx^j.
\end{array}\right.
\end{equation}

The main result of the generalized metrical multi-time Lagrange geometry is the theorem
of existence of the {\it Cartan canonical $h$-normal linear connection} $C\Gamma$ which  allow the
subsequent development of the  {\it generalized metrical  multi-time Lagrangian theory of
physical fields}.
\begin{th}{(of existence of Cartan canonical connection)}\\
On the generalized metrical multi-time Lagrange space $GML^n_p$ endowed with
the nonlinear connection $\Gamma$, there is a unique $h$-normal $\Gamma$-linear
connection, defined by adapted components
$$
C\Gamma=(H^\gamma_{\alpha\beta},G^k_{j\gamma},L^i_{jk},C^{i(\gamma)}
_{j(k)})
$$
having the metrical properties\medskip

i) $g_{ij\vert k}=0,\quad g_{ij}\vert^{(\gamma)}_{(k)}=0$,\medskip

ii) $\displaystyle{ G^k_{j\gamma}={g^{ki}\over 2}{\delta g_{ij}\over\delta
t^\gamma},\quad L^k_{ij}=L^k_{ji},\quad
C^{i(\gamma)}_{j(k)}=C^{i(\gamma)}_{k(j)}}$,\\
where $"_{/\beta}"$, $"_{\vert k}"$ and $"\vert_{(k)}^{(\gamma)}"$ are the local
covariant derivatives \cite{14} induced by $C\Gamma$.

Moreover, the coefficients $L^i_{jk}$ and $C^{i(\gamma)}_{j(k)}$ of the
Cartan canonical connection have the expressions
\begin{equation}
\begin{array}{l}
\displaystyle{L^i_{jk}={g^{im}\over 2}\left(
{\delta g_{mj}\over\delta x^k}+
{\delta g_{mk}\over\delta x^j}-
{\delta g_{jk}\over\delta x^m}\right),}\\
\displaystyle{C^{i(\gamma)}_{j(k)}={g^{im}\over 2}\left(
{\partial g_{mj}\over\partial x^k_\gamma}+
{\partial g_{mk}\over\partial x^j_\gamma}-
{\partial g_{jk}\over\partial x^m_\gamma}\right).}
\end{array}
\end{equation}
\end{th}
{\bf Proof.} Let $C\Gamma=(\bar G^\gamma_{\alpha\beta},G^k_{j\gamma},L^i_{jk},
C^{i(\gamma)}_{j(k)})$ be h-normal $\Gamma$-linear connection whose
coefficients are defined by
$\displaystyle{\bar G^\gamma_{\alpha\beta}=H^\gamma_{\alpha\beta},\;
G^k_{j\gamma}={g^{ki}\over 2}{\delta g_{ij}\over\delta t^\gamma},}$ and
\begin{equation}
\begin{array}{l}\medskip
\displaystyle{L^i_{jk}={g^{im}\over 2}\left({\delta g_{jm}\over\delta x^k}+
{\delta g_{km}\over\delta x^j}-{\delta g_{jk}\over\delta x^m}\right),}\\
\displaystyle{C^{i(\gamma)}_{j(k)}={g^{im}\over 2}\left({\partial g_{jm}\over\partial
x^k_\gamma}+{\partial g_{km}\over\partial x^j_\gamma}-{\partial g_{jk}\over
\partial x^m_\gamma}\right)}.
\end{array}
\end{equation}
By computations, one easily verifies that  $C\Gamma$ satisfies the  conditions
{\it i} and {\it  ii}.

Conversely, let  us consider  a h-normal
$\Gamma$-linear connection
$$
\tilde C\Gamma=(\tilde{\bar G}^\gamma_{\alpha\beta},
\tilde G^k_{j\gamma},\tilde L^i_{jk},\tilde C^{i(\gamma)}_{j(k)})
$$
which satisfies {\it  i} and {\it ii}. It follows
$$
\displaystyle{\tilde{\bar G}^\gamma_{\alpha\beta}=H^\gamma_
{\alpha\beta},\;\mbox{and}\;\;\;\tilde G^k_{j\gamma}={g^{ki}\over 2}{\delta g_{ij}\over
\delta t^\gamma}}.
$$

The condition $g_{ij\vert k}=0$ is equivalent to
$$
{\delta g_{ij}\over\delta x^k}=g_{mj}\tilde L^m_{ik}+g_{im}\tilde L^m_{jk}.
$$
Applying a Christoffel process to the indices $\{i,j,k\}$, we find
$$
\displaystyle{\tilde L^i_{jk}={g^{im}\over 2}\left({\delta g_{jm}\over\delta x^k}+
{\delta g_{km}\over\delta x^j}-{\delta g_{jk}\over\delta x^m}\right)}.
$$

By  analogy, using the relations $C^{i(\gamma)}_{j(k)}=C^{i(\gamma)}_{k(j)}$
and $g_{ij}\vert^{(\gamma)}_{(k)}=0$, following a Christoffel process applied
to the indices $\{i,j,k\}$, we obtain
$$
\displaystyle{\tilde C^{i(\gamma)}_{j(k)}={g^{im}\over 2}\left({\partial g_{jm}\over\partial
x^k_\gamma}+{\partial g_{km}\over\partial x^j_\gamma}-{\partial g_{jk}\over
\partial x^m_\gamma}\right)}.
$$

In conclusion, the uniqueness of the Cartan canonical  connection $C\Gamma$ is
clear. q. e. d.\medskip\\
\addtocounter{rem}{1}
{\bf Remark \therem}  As a rule, the  Cartan canonical connection of the generalized metrical
multi-time Lagrange space $GML^n_p$ verifies also the metrical properties
$$
h_{\alpha\beta/\gamma}=h_{\alpha\beta\vert k}=h_{\alpha\beta}\vert^{(\gamma)}
_{(k)}=0\;\mbox{and}\;g_{ij/\gamma}=0.
$$

Using the general results from \cite{14} upon the components  of the d-torsion
{\bf  T} and the d-curvature {\bf R} of an $h$-normal $\Gamma$-linear connection
$\nabla$ on $E$, we obtain without difficulties  the following theorems:

\begin{th}
The torsion d-tensor {\bf T} of the Cartan canonical connection of the generalized
metrical multi-time Lagrange space $GML^n_p$ is determined by eight local components
\begin{equation}
\begin{tabular}{|c|c|c|c|}
\hline
&$h_T$&$h_M$&$v$\\
\hline
$h_Th_T$&0&0&$R^{(m)}_{(\mu)\alpha\beta}$\\
\hline
$h_Mh_T$&0&$T^m_{\alpha j}$&$R^{(m)}_{(\mu)\alpha j}$\\
\hline
$h_Mh_M$&0&$0$&$R^{(m)}_{(\mu)ij}$\\
\hline
$vh_T$&0&0&$P^{(m)\;\;(\beta)}_{(\mu)\alpha(j)}$\\
\hline
$vh_M$&0&$P^{m(\beta)}_{i(j)}$&$P^{(m)\;(\beta)}_{(\mu)i(j)}$\\
\hline
$vv$&0&0&$S^{(m)(\alpha)(\beta)}_{(\mu)(i)(j)}$\\
\hline
\end{tabular}
\end{equation}
where\medskip

$
\displaystyle{P^{(m)\;\;(\beta)}_{(\mu)\alpha(j)}={\partial M^{(m)}_{(\mu)
\alpha}\over\partial x^j_\beta}-\delta^\beta_\mu G^m_{j\alpha}+\delta^m_j
H^\beta_{\mu\alpha},}\;\;
$\medskip

$\displaystyle{P^{(m)\;\;(\beta)}_{(\mu)i(j)}=
{\partial N^{(m)}_{(\mu)i}\over\partial x^j_\beta}-\delta^\beta_\mu L^m_{ji},}
$\medskip

$
\displaystyle{R^{(m)}_{(\mu)\alpha\beta}={\delta M^{(m)}_{(\mu)\alpha}\over
\delta t^\beta}-{\delta M^{(m)}_{(\mu)\beta}\over\delta t^\alpha},}
$\medskip

$\displaystyle{
R^{(m)}_{(\mu)\alpha j}={\delta M^{(m)}_{(\mu)\alpha}\over
\delta x^j}-{\delta N^{(m)}_{(\mu)j}\over\delta t^\alpha},}
$\medskip

$\displaystyle{
R^{(m)}_{(\mu)ij}={\delta N^{(m)}_{(\mu)i}\over\delta x^j}-{\delta N^{(m)}_
{(\mu)j}\over\delta x^i},}
$\medskip

$
S^{(m)(\alpha)(\beta)}_{(\mu)(i)(j)}=\delta^\alpha_\mu C^{m(\beta)}_{i(j)}-
\delta^\beta_\mu C^{m(\alpha)}_{j(i)}
$,\quad
$
T^m_{\alpha j}=-G^m_{j\alpha}$,\quad
$P^{m(\beta)}_{i(j)}=C^{m(\beta)}_{i(j)}$;
\end{th}

\begin{th}
The  curvature d-tensor {\bf R} of the Cartan
canonical connection of $GML^n_p$ is characterized by seven effective local
d-tensors,
\begin{equation}
\begin{tabular}{|c|c|c|c|}
\hline
&$h_T$&$h_M$&$v$\\
\hline
$h_Th_T$&$H^\alpha_{\eta\beta\gamma}$&$R^l_{i\beta\gamma}$&
$R^{(l)(\alpha)}_{(\eta)(i)\beta\gamma}=\delta^\alpha_\eta R^l_{i\beta\gamma}
+\delta^l_iH^\alpha_{\eta\beta\gamma}$\\
\hline
$h_Mh_T$&0&$R^l_{i\beta k}$&$R^{(l)(\alpha)}_{(\eta)(i)\beta k}=\delta^\alpha
_\eta R^l_{i\beta k}$\\
\hline
$h_Mh_M$&0&$R^l_{ijk}$&$R^{(l)(\alpha)}_{(\eta)(i)jk}=\delta^\alpha_\eta R^l_{
ijk}$\\
\hline
$vh_T$&0&$P^{l\;\;(\gamma)}_{i\beta(k)}$&$P^{(l)(\alpha)\;\;(\gamma)}_{(\eta)
(i)\beta(k)}=\delta^\alpha_\eta P^{l\;\;(\gamma)}_{i\beta(k)}$\\
\hline
$vh_M$&0&$P^{l\;(\gamma)}_{ij(k)}$&$P^{(l)(\alpha)\;(\gamma)}_{(\eta)
(i)j(k)}=\delta^\alpha_\eta P^{l\;(\gamma)}_{ij(k)}$\\
\hline
$vv$&0&$S^{l(\beta)(\gamma)}_{i(j)(k)}$&$S^{(l)(\alpha)(\beta)(\gamma)}_
{(\eta)(i)(j)(k)}=\delta^\alpha_\eta S^{l(\beta)(\gamma)}_{i(j)(k)}$\\
\hline
\end{tabular}
\end{equation}
where
$
\displaystyle{H^\alpha_{\eta\beta\gamma}={\partial H^\alpha_{\eta\beta}\over
\partial t^\gamma}-{\partial H^\alpha_{\eta\gamma}\over\partial t^\beta}+
H^\mu_{\eta\beta}H^\alpha_{\mu\gamma}-H^\mu_{\eta\gamma}H^\alpha_{\mu\beta},}
$\medskip

$
\displaystyle{R^l_{i\beta\gamma}={\delta G^l_{i\beta}\over\delta t^\gamma}-
{\delta G^l_{i\gamma}\over\delta t^\beta}+G^m_{i\beta}G^l_{m\gamma}-
G^m_{i\gamma}G^l_{m\beta}+C^{l(\mu)}_{i(m)}R^{(m)}_{(\mu)\beta\gamma},}
$\medskip

$
\displaystyle{R^l_{i\beta k}={\delta G^l_{i\beta}\over\delta x^k}-
{\delta L^l_{ik}\over\delta t^\beta}+G^m_{i\beta}L^l_{mk}-
L^m_{ik}G^l_{m\beta}+C^{l(\mu)}_{i(m)}R^{(m)}_{(\mu)\beta k},}
$\medskip

$
\displaystyle{R^l_{ijk}={\delta L^l_{ij}\over\delta x^k}-
{\delta L^l_{ik}\over\delta x^j}+L^m_{ij}L^l_{mk}-L^m_{ik}L^l_{mj}+
C^{l(\mu)}_{i(m)}R^{(m)}_{(\mu)jk},}
$\medskip

$
\displaystyle{P^{l\;\;(\gamma)}_{i\beta(k)}={\partial G^l_{i\beta}\over\partial
x^k_\gamma}-C^{l(\gamma)}_{i(k)/\beta}+C^{l(\mu)}_{i(m)}P^{(m)\;\;(\gamma)}_
{(\mu)\beta(k)},}
$\medskip

$
\displaystyle{P^{l\;(\gamma)}_{ij(k)}={\partial L^l_{ij}\over\partial
x^k_\gamma}-C^{l(\gamma)}_{i(k)\vert j}+C^{l(\mu)}_{i(m)}P^{(m)\;(\gamma)}_
{(\mu)j(k)},}
$\medskip

$
\displaystyle{S^{l(\beta)(\gamma)}_{i(j)(k)}={\partial C^{l(\beta)}_{i(j)}
\over\partial x^k_\gamma}-{\partial C^{l(\gamma)}_{i(k)}\over\partial x^j_
\beta}+C^{m(\beta)}_{i(j)}C^{l(\gamma)}_{m(k)}-
C^{m(\gamma)}_{i(k)}C^{l(\beta)}_{m(j)}.}
$
\end{th}

Moreover, in the particular case of Cartan canonical connection, applying the
Ricci identities attached to an {\it  $h$-normal $\Gamma$-linear connection
of Cartan type} for an arbitrary d-tensor \cite{11} to the  metric
d-tensors $h_{\alpha\beta}(t^\gamma)$ and $g_{ij}(t^\gamma,x^k,x^k_\gamma)$
of $GML^n_p$, it follows
\begin{th}
The following curvature d-tensor identities are true:
\begin{equation}\label{ci}
\begin{array}{lll}\medskip
1)&H_{\alpha\beta\gamma\delta}+H_{\beta\alpha\gamma\delta}=0,&H_{\alpha\beta
\gamma\delta}=h_{\beta\mu}H^\mu_{\alpha\gamma\delta};\\\medskip
2)&R_{ij\beta\gamma}+R_{ji\beta\gamma}=0,&R_{ij\beta\gamma}=g_{jm}R^m_{i\beta
\gamma};\\\medskip
3)&R_{ij\beta k}+R_{ji\beta k}=0,&R_{ij\beta k}=g_{jm}R^m_{i\beta k};
\\\medskip
4)&R_{ijkl}+R_{jikl}=0,&R_{ijkl}=g_{jm}R^m_{ikl};
\\\medskip
5)&P_{ij\beta(k)}^{\;\;\;\;(\gamma)}+P_{ji\beta(k)}^{\;\;\;\;(\gamma)}=0,&
P_{ij\beta(k)}^{\;\;\;\;(\gamma)}=g_{jm}P^{m\;(\gamma)}_{i\beta(k)};
\\\medskip
6)&P_{ijk(l)}^{\;\;\;\;(\delta)}+P_{jik(l)}^{\;\;\;\;(\delta)}=0,&P_{ijk(l)}^
{\;\;\;\;(\delta)}=g_{jm}P^{m\;(\delta)}_{ik(l)};
\\\medskip
7)&S_{ij(k)(l)}^{\;\;(\gamma)(\delta)}+S_{ji(k)(l)}^{\;\;(\gamma)(\delta)}=0,&
S_{ij(k)(l)}^{\;\;(\gamma)(\delta)}=g_{jm}S^{m(\gamma)(\delta)}_{i(k)(l)}.
\end{array}
\end{equation}
\end{th}

\section{Electromagnetic field. Maxwell equations}

\setcounter{equation}{0}
\hspace{5mm} Let $GML^n_p=(J^1(T,M),G^{(\alpha)(\beta)}_{(i)(j)}(t^\gamma,x^k,
x^k_\gamma)=h^{\alpha\beta}(t^\gamma)g_{ij}(t^\gamma,x^k,x^k_\gamma))$ be a
generalized metrical multi-time Lagrange  space and let us consider $M^{(i)}_{(\alpha)\beta}=
-H^\mu_{\alpha\beta}x^i_\mu$ the canonical temporal nonlinear connection of $GML^n_p$.
Let us suppose  that $GML^n_p$ is endowed  with an {\it without torsion}
spatial nonlinear connection $N^{(i)}_{(\alpha)j}$, that  is, it is verified
the relation
\begin{equation}
\displaystyle{{\partial N^{(i)}_{(\alpha)j}\over\partial x^k_\gamma}=
{\partial N^{(i)}_{(\alpha)k}\over\partial x^j_\gamma}}.
\end{equation}
\addtocounter{ex}{1}
{\bf Examples \theex} i) It is easy to verify that the spatial nonlinear
connection
$$
N^{(i)}_{(\alpha)j}=\gamma^i_{jm}x^m_\alpha,
$$
used in the examples {\bf 2.5} and {\bf 2.6}, is without torsion.

ii) The spatial nonlinear connection
$$
N^{(i)}_{(\alpha)j}=\Gamma^i_{jm}x^m_\alpha+{g^{im}\over 2}{\partial g_{jm}
\over\partial t^\alpha},
$$
from the example {\bf 2.3} is also without torsion.\medskip\\
Let us denote $C\Gamma=(H^\gamma_{\alpha\beta},G^k_{i\gamma},L^k_{ij},
C^{k(\gamma)}_{i(j)})$ the Cartan canonical connection of $GML^n_p$, and
$"_{/\beta}"$, $"_{\vert k}"$, $"\vert^{(\gamma)}_{(k)}"$ its local covariant
derivatives.\medskip\\
\addtocounter{rem}{1}
{\bf Remark \therem} Because we use an without torsion spatial  nonlinear
connection, the torsion d-tensor
$$
\displaystyle{P^{(m)\;\;(\beta)}_{(\mu)i(j)}=
{\partial N^{(m)}_{(\mu)i}\over\partial x^j_\beta}-\delta^\beta_\mu L^m_{ji}}
$$
of the Cartan canonical connection is symmetric in $i$ and $j$. In other
words, we have $P^{(m)\;\;(\beta)}_{(\mu)i(j)}=P^{(m)\;\;(\beta)}_{(\mu)j(i)}$.
\medskip

In order to develope the generalized metrical multi-time theory of
electromagnetism, we construct the {\it metrical deflection d-tensors}
\begin{equation}
\begin{array}{l}\medskip
\bar D^{(\alpha)}_{(i)\beta}=\left[G^{(\alpha)(\mu)}_{(i)(m)}x^m_\mu\right]_
{/\beta},\\\medskip
D^{(\alpha)}_{(i)j}=\left[G^{(\alpha)(\mu)}_{(i)(m)}x^m_\mu\right]_{\vert j},
\\\medskip
d^{(\alpha)(\beta)}_{(i)(j)}=\left[G^{(\alpha)(\mu)}_{(i)(m)}x^m_\mu\right]
\vert^{(\beta)}_{(j)},
\end{array}
\end{equation}
where $G^{(\alpha)(\beta)}_{(i)(k)}$ is the vertical fundamental  metrical
d-tensor  of the generalized metrical multi-time Lagrange space $GML^n_p$,
and $x^m_\mu$ are  the components of the {\it canonical Liouville d-tensor}
field {\bf C}$\displaystyle{=x^m_\mu{\partial\over\partial x^m_\mu}}$.
\medskip\\
\addtocounter{defin}{1}
{\bf Definition \thedefin} The distinguished 2-form on $J^1(T,M)$,
\begin{equation}
F=F^{(\alpha)}_{(i)j}\delta x^i_\alpha\wedge dx^j+f^{(\alpha)(\beta)}_{(i)(j)}
\delta x^i_\alpha\wedge\delta x^j_\beta,
\end{equation}
where
$F^{(\alpha)}_{(i)j}=\displaystyle{{1\over 2}\left[D^{(\alpha)}
_{(i)j}-D^{(\alpha)}_{(j)i}\right]}$ and
$f^{(\alpha)(\beta)}_{(i)(j)}=\displaystyle{{1\over 2}\left[d^{(\alpha)(\beta)}
_{(i)(j)}-d^{(\alpha)(\beta)}_{(j)(i)}\right]}$, is called the distinguished
{\it electromagnetic 2-form} of the generalized  metrical multi-time Lagrange
space $GML^n_p$.\medskip\\
\addtocounter{rem}{1}
{\bf Remark \therem} The previous definition is  a  natural generalization
of the electromagnetic form of the Miron-Anastasiei electromagnetism
\cite{7}.\medskip

The main result of the generalized metrical  multi-time Lagrangian
theory of electromagnetism is the following
\begin{th}
The electromagnetic components $F^{(\alpha)}_{(i)j}$ and $f^{(\alpha)(\beta)}_
{(i)(j)}$ of the generalized metrical multi-time Lagrange space $GML^n_p$ are
governed by the Maxwell equations:
\medskip
$$
\left\{\begin{array}{l}\medskip
\displaystyle{F^{(\alpha)}_{(i)k/\beta}={1\over 2}{\cal A}_
{\{i,k\}}\left\{\bar D^{(\alpha)}_{(i)\beta\vert k}
+D^{(\alpha)}_{(i)m}T^m_{\beta k}+d^{(\alpha)(\mu)}_
{(i)(m)}R^{(m)}_{(\mu)\beta k}-\right.}\\\medskip
\displaystyle{\hspace*{30mm}-\left.\left[T^p_{\beta i\vert k}+C^{p(\mu)}_
{k(m)}R^{(m)}_{(\mu)\beta i}\right]
x^{(\alpha)}_{(p)}\right\}}\\\medskip
\displaystyle{f^{(\alpha)(\gamma)}_{(i)(k)/\beta}={1\over 2}{\cal A}_{\{i,k\}}
\left\{\bar D^{(\alpha)}_{(i)\beta}\vert^{(\gamma)}_{(k)}+d^{(\alpha)(\mu)}_
{(i)(m)}P^{(m)\;(\gamma)}_{(\mu)\beta(k)}-\right.}\\\medskip
\displaystyle{\hspace*{30mm}-\left.\left[{\partial T^p_{\beta i}\over
\partial x^k_\gamma}+C^{p(\mu)}_{k(m)}P^{(m)\;(\gamma)}_{(\mu)\beta(i)}\right]
x^{(\alpha)}_{(p)}\right\}}\\\medskip
\displaystyle{\sum_{\{i,j,k\}}F^{(\alpha)}_{(i)j\vert k}=-{1\over 2}\sum_
{\{i,j,k\}}}\left[C^{p(\mu)}_{i(m)}x^{(\alpha)}_{(p)}+d^{(\alpha)(\mu)}_{(i)(m)}
\right]R^{(m)}_{(\mu)jk}\\\medskip
\displaystyle{\sum_{\{i,j,k\}}\left\{F^{(\alpha)}_{(i)j}\vert^{(\gamma)}_{(k)}+
f^{(\alpha)(\gamma)}_{(i)(j)\vert k}\right\}=0}
\\
\displaystyle{\sum_{\{i,j,k\}}f^{(\alpha)(\beta)}_{(i)(j)}\vert^{(\gamma)}_
{(k)}=0,}
\end{array}\right.
$$
where $x^{(\alpha)}_{(p)}=G^{(\alpha)(\mu)}_{(p)(m)}x^m_\mu$.
\end{th}
{\bf Proof.} Now, let us consider the following general deflection d-tensor
identities \cite{14}\medskip

$$
\begin{array}{ll}\medskip
d_1)&\bar D^{(p)}_{(\nu)\beta\vert k}-D^{(p)}_{(\nu)k/\beta}=x^m_\nu
R^p_{m\beta k}-D^{(p)}_{(\nu)m}T^m_{\beta k}-d^{(p)(\mu)}_{(\nu)(m)}
R^{(m)}_{(\mu)\beta k},\\\medskip
d_2)&\bar D^{(p)}_{(\nu)\beta}\vert^{(\gamma)}_{(k)}-d^{(p)(\gamma)}_{(\nu)
(k)/\beta}=x^m_\nu P^{p\;\;(\gamma)}_{m\beta(k)}-d^{(p)(\mu)}
_{(\nu)(m)}P^{(m)\;\;(\gamma)}_{(\mu)\beta(k)},\\\medskip
d_3)&D^{(p)}_{(\nu)j\vert k}-D^{(p)}_{(\nu)k\vert j}=x^m_\nu R^p_{mjk}-
d^{(p)(\mu)}_{(\nu)(m)}R^{(m)}_{(\mu)jk},
\\\medskip
d_4)&D^{(p)}_{(\nu)j}\vert^{(\gamma)}_{(k)}-d^{(p)(\gamma)}_{(\nu)
(k)\vert j}=x^m_\nu P^{p\;(\gamma)}_{mj(k)}-D^{(p)}_{(\nu)m}
C^{m(\gamma)}_{j(k)}-d^{(p)(\mu)}_{(\nu)(m)}P^{(m)\;(\gamma)}_
{(\mu)j(k)},\\\medskip
d_5)&d^{(p)(\beta)}_{(\nu)(j)}\vert^{(\gamma)}_{(k)}-d^{(p)(\gamma)}_{(\nu)
(k)}\vert^{(\beta)}_{(j)}=x^m_\nu S^{p(\beta)(\gamma)}_{m(j)(k)}-
d^{(p)(\mu)}_{(\nu)(m)}S^{(m)(\beta)(\gamma)}_{(\mu)(j)(k)},
\end{array}
$$
where
$\bar D^{(i)}_{(\alpha)\beta}=x^i_{\alpha/\beta}$,
$D^{(i)}_{(\alpha)j}=x^i_{\alpha\vert j}$,
$d^{(i)(\beta)}_{(\alpha)(j)}=x^i_\alpha\vert^{(\beta)}_{(j)}$.
Contracting the above deflection d-tensor identities by $G^{(\alpha)(\nu)}_{(i)(p)}$
and using the curvature d-tensor identities  \ref{ci}, we obtain the following
metrical deflection d-tensors  identities:\medskip
$$
\begin{array}{ll}\medskip
d^\prime_1)&\bar D^{(\alpha)}_{(i)\beta\vert k}-D^{(\alpha)}_{(i)k/\beta}=
-x^{(\alpha)}_{(m)}R^m_{i\beta k}-D^{(\alpha)}_{(i)m}T^m_{\beta k}-d^{(\alpha)
(\mu)}_{(i)(m)}R^{(m)}_{(\mu)\beta k},\\\medskip
d^\prime_2)&\bar D^{(\alpha)}_{(i)\beta}\vert^{(\gamma)}_{(k)}-d^{(\alpha)
(\gamma)}_{(i)(k)/\beta}=-x^{(\alpha)}_{(m)}P^{m\;\;(\gamma)}_{i\beta(k)}-
d^{(\alpha)(\mu)}_{(i)(m)}P^{(m)\;\;(\gamma)}_{(\mu)\beta(k)},\\\medskip
d^\prime_3)&D^{(\alpha)}_{(i)j\vert k}-D^{(\alpha)}_{(i)k\vert j}=-x^{(\alpha)}
_{(m)}R^m_{ijk}-d^{(\alpha)(\mu)}_{(i)(m)}R^{(m)}_{(\mu)jk},
\\\medskip
d^\prime_4)&D^{(\alpha)}_{(i)j}\vert^{(\gamma)}_{(k)}-d^{(\alpha)(\gamma)}
_{(i)(k)\vert j}=-x^{(\alpha)}_{(m)}P^{m\;(\gamma)}_{ij(k)}-D^{(\alpha)}_{(i)
m}C^{m(\gamma)}_{j(k)}-d^{(\alpha)(\mu)}_{(i)(m)}P^{(m)\;(\gamma)}_{(\mu)j(k)},
\\\medskip
d^\prime_5)&d^{(\alpha)(\beta)}_{(i)(j)}\vert^{(\gamma)}_{(k)}-d^{(\alpha)
(\gamma)}_{(i)(k)}\vert^{(\beta)}_{(j)}=-x^{(\alpha)}_{(m)}S^{m(\beta)(\gamma)}
_{i(j)(k)}-d^{(\alpha)(\mu)}_{(i)(m)}S^{(m)(\beta)(\gamma)}_{(\mu)(j)(k)}.
\end{array}
$$

At the same time, we recall that the following Bianchi identities \cite{11}\medskip
$$
\begin{array}{ll}\medskip
b_1)&{\cal A}_{\{j,k\}}\left\{R^l_{j\alpha k}+T^l_{\alpha j\vert k}+
C^{l(\mu)}_{k(m)}R^{(m)}_{(\mu)\alpha j}\right\}=0,\medskip\\
b_2)&T^l_{\alpha k}\vert^{(\varepsilon)}_{(p)}-C^{l(\varepsilon)}_{m(p)}
T^m_{\alpha k}+P^{l\;\;(\varepsilon)}_{k\alpha(p)}-C^{l(\varepsilon)}_
{k(p)/\alpha}-C^{l(\mu)}_{k(m)}P^{(m)\;\;(\varepsilon)}_{(\mu)\alpha(p)}=0,
\medskip\\
b_3)&\sum_{\{i,j,k\}}\left\{R^l_{ijk}-C^{l(\mu)}_{k(m)}R^{(m)}_{(\mu)ij}
\right\}=0,\medskip\\
b_4)&{\cal A}_{\{j,k\}}\left\{P^{l\;\;(\varepsilon)}_{jk(p)}+C^{l(\varepsilon)}
_{j(p)\vert k}+C^{l(\mu)}_{k(m)}P^{(m)\;\;(\varepsilon)}_{(\mu)j(p)}\right\}=0,
\\
b_5)&\displaystyle{S^{l(\beta)(\gamma)}_{i(j)(k)}={\partial C^{l(\beta)}_{i(j)}
\over\partial x^k_\gamma}-{\partial C^{l(\gamma)}_{i(k)}\over\partial x^j_
\beta}+C^{m(\beta)}_{i(j)}C^{l(\gamma)}_{m(k)}-
C^{m(\gamma)}_{i(k)}C^{l(\beta)}_{m(j)}},
\end{array}
$$
where ${\cal A}_{\{j,k\}}$ means alternate sum and $\sum_{\{i,j,k\}}$ means
cyclic sum, hold good.

In order to obtain the first Maxwell identity, we permute $i$ and $k$ in
$d^\prime_1$ and we subtract the new identity from the initial one. Finally,
using the Bianchi identity $b_1$, we obtain what we were looking for. By
analogy, using  $d_2^\prime$  and $b_2$, it follows the  second  Maxwell
equation.

Doing a cyclic sum by the indices $\{i,j,k\}$ in $d^\prime_3$ and using the
Bianchi identity $b_3$, it follows the third Maxwell equation.

Applying a Christoffel process to the indices $\{i,j,k\}$ in $d^\prime_4$
and combining with the Bianchi identity $b_4$ and the relation
$P^{(m)\;\;(\varepsilon)}_{(\mu)j(p)}=P^{(m)\;\;(\varepsilon)}_{(\mu)p(j)}$,
we get a new identity. The cyclic sum by the indices
$\{i,j,k\}$ applied to this last identity implies the fourth Maxwell equation.

Using $b_5$ and the relation  $S^{(m)(\alpha)(\beta)}_{(\mu)(i)(j)}=\delta^\alpha
_\mu C^{m(\beta)}_{i(j)}-\delta^\mu_\beta C^{m(\alpha)}_{j(i)}$, a Christoffel
process applied in $d^\prime_5$ gives a new identity. Doing a cyclic sum by
the indices $\{i,j,k\}$ in this identity, we obtain the last Maxwell equation.
\linebreak
q. e. d.\medskip\\
\addtocounter{rem}{1}
{\bf Remark \therem} Let us suppose that $\dim T\geq 2$. In the particular
case of a $GML^n_p$ with $g_{ij}(t^\gamma,x^k,x^k_\gamma)=g_{ij}(t^\gamma,x^k)$,
the Maxwell equations reduce to that of a metrical multi-time Lagrange space
\cite{13}.

\section{Gravitational field. Einstein equations}

\setcounter{equation}{0}
\hspace{5mm} In order to develope the generalized metrical multi-time
Lagrange theory of gravitational field, we introduce the following\medskip\\
\addtocounter{defin}{1}
{\bf Definition \thedefin} From physical point of view, an adapted metrical d-tensor
$G$ on $E=J^1(T,M)$, expressed locally by
$$
G=h_{\alpha\beta}dt^\alpha\otimes dt^\beta+g_{ij}dx^i\otimes dx^j+h^{\alpha
\beta}g_{ij}\delta x^i_\alpha\otimes\delta x^j_\beta,
$$
where $g_{ij}=g_{ij}(t^\gamma,x^k,x^k_\gamma)$ is a d-tensor field on $E$,
symmetric, of rank $n=\dim M$ and having a constant signature on $E$, is
called a {\it gravitational $h$-potential}.\medskip\\
\addtocounter{rem}{1}
{\bf Remark \therem} The naturalness of this definition comes from the particular
case $(T,h)=(R,\delta)$. In this case, we recover the {\it gravitational
potentials} $g_{ij}(x,y)$ from Miron-Anastasiei theory of gravitational field
\cite{7}.\medskip

Now, let $GML^n_p=(J^1(T,M),G^{(\alpha)(\beta)}_{(i)(j)}=h^{\alpha\beta}(t^
\gamma)g_{ij}(t^\gamma,x^k,x^k_\gamma))$ be a
generalized metrical multi-time Lagrange  space and let $M^{(i)}_{(\alpha)\beta}=
-H^\mu_{\alpha\beta}x^i_\mu$ be  the canonical temporal nonlinear connection of $GML^n_p$.
Let us suppose  that $GML^n_p$ is endowed  with a
spatial nonlinear connection $N^{(i)}_{(\alpha)j}$, not necessarily without
torsion.
It is obvious that the vertical fundamental metrical d-tensor of $GML^n_p$
induces a natural gravitational $h$-potential (i. e. a  Sasakian-like metric on
$J^1(T,M)$), setting
$$
G=h_{\alpha\beta}dt^\alpha\otimes dt^\beta+g_{ij}dx^i\otimes dx^j+h^{\alpha
\beta}g_{ij}\delta x^i_\alpha\otimes\delta x^j_\beta.
$$
Let us consider $C\Gamma=(H^\gamma_{\alpha\beta},G^k_{j\gamma},L^i_{jk},
C^{i(\gamma)}_{j(k)})$ the Cartan canonical connection of $GML^n_p$.

We postulate that the Einstein equations which govern the gravitational
$h$-potential $G$ of the generalized metrical multi-time Lagrange space
$GML^n_p$ are the Einstein equations attached to the Cartan canonical
connection $C\Gamma$ of $GML^n_p$ and the adapted metric $G$ on $E$, that is,
\begin{equation}
Ric(C\Gamma)-{Sc(C\Gamma)\over 2}G={\cal K}{\cal T},
\end{equation}
where $Ric(C\Gamma)$ represents the Ricci d-tensor of the Cartan connection,
$Sc(C\Gamma)$ is its scalar curvature, ${\cal K}$ is the Einstein constant and ${\cal T}$
is an intrinsec tensor of matter which is called  the {\it stress-energy}
d-tensor.

In the adapted basis $(X_A)=\displaystyle{\left({\delta\over\delta t^\alpha},
{\delta\over\delta x^i},{\partial\over\partial x^i_\alpha}\right)}$ of  the
nonlinear  connection $\Gamma$ of $GML^n_p$, the curvature d-tensor {\bf R} of
the Cartan connection is expressed locally by {\bf R}$(X_C,X_B)X_A=R^D_{ABC}X_D$.
It follows that we have $R_{AB}=Ric(X_A,X_B)=R^D_{ABD}$ and $Sc(C\Gamma)
=G^{AB}R_{AB}$, where
\begin{equation}
G^{AB}=\left\{\begin{array}{ll}\medskip
h_{\alpha\beta},&\mbox{for}\;\;A=\alpha,\;B=\beta\\\medskip
g^{ij},&\mbox{for}\;\;A=i,\;B=j\\\medskip
h_{\alpha\beta}g^{ij},&\mbox{for}\;\;A={(i)\atop(\alpha)},\;B={(j)\atop(\beta)}\\
0,&\mbox{otherwise}.
\end{array}\right.
\end{equation}

Taking into account, on the one hand, the form of the vertical fundamental
metrical  d-tensor  $G^{(\alpha)(\beta)}_{(i)(j)}$   of the generalized
metrical multi-time Lagrange space $GML^n_p$, and,  on the other  hand,
the expressions of local curvature d-tensors attached to the  Cartan
canonical  connection $C\Gamma$, by a direct calculation, we deduce
\begin{th}
The Ricci d-tensor $Ric(C\Gamma)$ of the  Cartan canonical connection $C\Gamma$
of the generalized metrical multi-time Lagrange space  $GML^n_p$, is determined
by the following components:
$$
\begin{array}{l}\medskip
R_{(\alpha)(\beta)}\stackrel{not}{=}H_{\alpha\beta}=H^\mu_{\alpha\beta\mu},\quad
R^{\;(\alpha)}_{i(j)}\stackrel{not}{=}P^{\;(\alpha)}_{i(j)}=-P^{m\;(\alpha)}
_{im(j)},\\\medskip
R^{(\alpha)}_{(i)j}\stackrel{not}{=}P^{(\alpha)}_{(i)j}=P^{m\;(\alpha)}_{ij(m)},
\quad
R^{(\alpha)}_{(i)\beta}\stackrel{not}{=}P^{(\alpha)}_{(i)\beta}=P^{m\;(\alpha)}_
{i\beta(m)},\\\medskip
R^{(\alpha)(\beta)}_{(i)(j)}\stackrel{not}{=}S^{(\alpha)(\beta)}_{(i)(j)}=
S^{m(\beta)(\alpha)}_{i(j)(m)},\quad
R_{i\alpha}=R^m_{i\alpha m},\quad R_{ij}=R^m_{ijm}.
\end{array}
$$
\end{th}
\begin{cor}
The scalar curvature  $Sc(C\Gamma)$ of the Cartan canonical connection $C\Gamma$
of the generalized metrical multi-time Lagrange space $GML^n_p$, is given  by
$$
Sc(C\Gamma)=H+R+S,
$$
where $H=h^{\alpha\beta}H_{\alpha\beta},\;R=g^{ij}R_{ij}$
and $S=h_{\alpha\beta}g^{ij}S^{(\alpha)\beta)}_{(i)(j)}$.
\end{cor}

The main result of the generalized metrical multi-time Lagrange theory of
gravitational field is  given by
\begin{th}
The Einstein equations which govern the gravitational $h$-potential
$G$ of the generalized metrical multi-time Lagrange space $GML^n_p$,
have the local form
$$
\left\{\begin{array}{l}\medskip
\displaystyle{H_{\alpha\beta}-{H+R+S\over 2}h_{\alpha\beta}={\cal K}{\cal T}_
{\alpha\beta}}\\\medskip
\displaystyle{R_{ij}-{H+R+S\over 2}g_{ij}={\cal K}{\cal T}_{ij}}\\\medskip
\displaystyle{S^{(\alpha)(\beta)}_{(i)(j)}-{H+R+S\over 2}h^{\alpha\beta}g_{ij}={\cal K}{\cal T}^{(\alpha)
(\beta)}_{(i)(j)}},
\end{array}\right.\leqno{(E_1)}
$$
$$
\left\{\begin{array}{lll}\medskip
0={\cal T}_{\alpha i},&R_{i\alpha}={\cal K}{\cal T}_{i\alpha},&
P^{(\alpha)}_{(i)\beta}={\cal K}{\cal T}^{(\alpha)}_{(i)\beta}\\
0={\cal T}^{\;(\beta)}_{\alpha(i)},&
P^{\;(\alpha)}_{i(j)}={\cal K}{\cal T}^{\;(\alpha)}_{i(j)},&
P^{(\alpha)}_{(i)j}={\cal K}{\cal T}^{(\alpha)}_{(i)j},
\end{array}\right.\leqno{(E_2)}
$$
where ${\cal T}_{AB},\;A,B\in\{\alpha,i,{(\alpha)\atop(i)}\}$ are the adapted
local components of the stress-energy d-tensor ${\cal T}$.
\end{th}
\addtocounter{rem}{1}
{\bf Remarks \therem} i) In order to have the compatibility of the Einstein equations, it is
necessary that the certain  adapted local components of the stress-energy
d-tensor vanish {\it "a priori"}.

ii)  If $p=\dim T\geq 2$ and $g_{ij}=g_{ij}(t^\gamma,x^k)$, the Einstein equations
of $GML^n_p$ reduce to the  Einstein equations of a metrical multi-time Lagrange
space \cite{13}.\medskip

From physical point of view, it is  well known that the stress-energy
d-tensor ${\cal T}$ must  verify the local {\it conservation laws} ${\cal T}^
B_{A\vert B}=0,\;\forall\;A\in\{\alpha,i,{(\alpha)\atop (i)}\}$,
where ${\cal T}^B_A=G^{BD}{\cal T}_{DA}$. Consequently, by a direct calculation,
we find the  following
\begin{th}
In  the generalized metrical multi-time Lagrange space $GML^n_p$,  the  following
conservation laws of the Einstein equations hold good:
\begin{equation}
\left\{\begin{array}{l}\medskip
\displaystyle{\left[H^\mu_\beta-{H+R+S\over 2}\delta^\mu_\beta\right]_{/\mu}=
-R^m_{\beta\vert m}-P^{(m)}_{(\mu)\beta}\vert^{(\mu)}_{(m)}}\\\medskip
\displaystyle{\left[R^m_j-{H+R+S\over 2}\delta^m_j\right]_{\vert m}=-P^{(m)}_
{(\mu)j}\vert^{(\mu)}_{(m)}}\\
\displaystyle{\left[S^{(m)(\beta)}_{(\mu)(j)}-{H+R+S\over 2}\delta^m_j\delta^\beta
_\mu\right]\vert^{(\mu)}_{(m)}=-P^{m(\beta)}_{\;\;\;\;j\vert m}},
\end{array}\right.
\end{equation}
where
\begin{equation}\label{dt}
\begin{array}{l}\medskip
H^\alpha_\beta=h^{\alpha\mu}H_{\mu\beta},\;\;R^i_\beta=g^{im}R_{m\beta},\;\;
R^i_j=g^{im}R_{mj},\\\medskip
P^{(i)}_{(\alpha)\beta}=g^{im}h_{\alpha\mu}P^{(\mu)}_{(m)\beta},\;\;
P^{i(\beta)}_{\;\;(j)}=g^{im}P_{m(j)}^{\;\;(\beta)},\\
P^{(i)}_{(\alpha)j}=g^{im}h_{\alpha\mu}P^{(\mu)}_{(m)j},\;\;
S^{(i)(\beta)}_{(\alpha)(j)}=g^{im}h_{\alpha\mu}S^{(\mu)(\beta)}_{(m)(j)}.
\end{array}
\end{equation}
\end{th}

\section{A natural form of Einstein equations}

\setcounter{equation}{0}
\hspace{5mm} Let us suppose that $p=\dim T>2$ and $n=\dim M>2$. In this context,
we will show that the  Einstein equations  and their  conservation laws can
be  rewritten in  a  more natural form.
\begin{th}
The Einstein equations $(E_1)$ of $GML^n_p$ are equivalent  to the  set of
equations
$$
\left\{\begin{array}{l}\medskip
\displaystyle{H_{\alpha\beta}-{H\over 2}h_{\alpha\beta}={\cal K}\tilde{\cal T}_
{\alpha\beta}}\\\medskip
\displaystyle{R_{ij}-{R\over 2}g_{ij}={\cal K}\tilde{\cal T}_{ij}}\\
\displaystyle{S^{(\alpha)(\beta)}_{(i)(j)}-{S\over 2}h^{\alpha\beta}g_{ij}=
{\cal K}\tilde{\cal T}^{(\alpha)(\beta)}_{(i)(j)}},
\end{array}\right.\leqno{(E^\prime_1)}
$$
where $\tilde{\cal T}_{\alpha\beta}$, $\tilde{\cal T}_{ij}$ and $\tilde{\cal
T}^{(\alpha)(\beta)}_{(i)(j)}$  represent the components  of new stress-energy
d-tensor ${\cal T}$.
\end{th}
{\bf Proof.} $(E_1)\Longrightarrow (E^\prime_1)$ Contracting the  equations
$(E_1)$, in order by $h^{\alpha\beta}$, $g^{ij}$ and $G^{(\alpha)(\beta)}_{(i)(j)}$, we
obtain the  system
\begin{equation}\label{sys}
\left\{\begin{array}{l}\medskip
\displaystyle{H-{p\over 2}(H+R+S)={\cal K}{\cal T}_T}\\\medskip
\displaystyle{R-{n\over 2}(H+R+S)={\cal K}{\cal T}_M}\\
\displaystyle{S-{pn\over 2}(H+R+S)={\cal K}{\cal T}_v},
\end{array}\right.
\end{equation}
where ${\cal T}_T=h^{\alpha\beta}{\cal T}_{\alpha\beta}$, ${\cal T}_M=g^{ij}
{\cal T}_{ij}$ and ${\cal T}_v=G^{(m)(r)}_{(\mu)(\nu)}{\cal T}^{(\mu)(\nu)}_
{(m)(r)}$. Solving the algebraic system \ref{sys} in the unknowns $H$, $R$ and
$S$, we find
\begin{equation}
\left\{\begin{array}{l}\medskip
\displaystyle{H={\cal K}\left[{\cal  T}_T+{p\over 2-p-n-pn}({\cal T}_T+{\cal T}
_M+{\cal T}_v)\right]}\\\medskip
\displaystyle{R={\cal K}\left[{\cal  T}_M+{n\over 2-p-n-pn}({\cal T}_T+{\cal T}
_M+{\cal T}_v)\right]}\\
\displaystyle{S={\cal K}\left[{\cal  T}_v+{pn\over 2-p-n-pn}({\cal T}_T+{\cal T}
_M+{\cal T}_v)\right]}.
\end{array}\right.
\end{equation}
Defining the components of the new stress-energy d-tensor ${\cal T}$ by
\begin{equation}\label{*}
\left\{\begin{array}{l}\medskip
\displaystyle{\tilde{\cal T}_{\alpha\beta}={\cal T}_{\alpha\beta}+{R+S\over
2{\cal K}}h_{\alpha\beta}}\\\medskip
\displaystyle{\tilde{\cal T}_{ij}={\cal T}_{ij}+{H+S\over 2{\cal K}}g_{ij}}\\
\displaystyle{\tilde{\cal T}^{(\alpha)(\beta)}_{(i)(j)}={\cal T}^{(\alpha)
(\beta)}_{(i)(j)}+{H+R\over 2{\cal K}}G^{(\alpha)(\beta)}_{(i)(j)}},
\end{array}\right.
\end{equation}
we obtain the new  form $(E_1^\prime)$  of the Einstein equations $(E_1)$.\medskip

$(E^\prime_1)\Longrightarrow (E_1)$ Again by contractions,  the  system $(E_1^\prime)$
implies the equalities
\begin{equation}\label{**}
\left\{\begin{array}{l}\medskip
\displaystyle{H={2{\cal K}\tilde{\cal T}_T\over 2-p}}\\\medskip
\displaystyle{R={2{\cal K}\tilde{\cal T}_M\over 2-n}}\\
\displaystyle{S={2{\cal K}\tilde{\cal T}_v\over 2-pn}},
\end{array}\right.
\end{equation}
where $\tilde{\cal T}_T=h^{\alpha\beta}\tilde{\cal T}_{\alpha\beta}$, $\tilde
{\cal T}_M=g^{ij}\tilde{\cal T}_{ij}$ and $\tilde{\cal T}_v=G^{(m)(r)}_{(\mu)
(\nu)}\tilde{\cal T}^{(\mu)(\nu)}_{(m)(r)}$. In the  sequel, setting  the
components  of the  stress-energy d-tensor  ${\cal T}$ as
\begin{equation}
\left\{\begin{array}{l}\medskip
\displaystyle{{\cal T}_{\alpha\beta}=\tilde{\cal T}_{\alpha\beta}-{R+S\over
2{\cal K}}h_{\alpha\beta}}\\\medskip
\displaystyle{{\cal T}_{ij}=\tilde{\cal T}_{ij}-{H+S\over 2{\cal K}}g_{ij}}\\
\displaystyle{{\cal T}^{(\alpha)(\beta)}_{(i)(j)}=\tilde{\cal T}^{(\alpha)
(\beta)}_{(i)(j)}-{H+R\over 2{\cal K}}G^{(\alpha)(\beta)}_{(i)(j)}},
\end{array}\right.
\end{equation}
we obtain  what we  were looking for. q. e. d.\medskip

Now, let us  study the form of conservation laws associated to the new
stress-energy d-tensor $\tilde{\cal T}$. In order to  describe these new conservation
laws, we use  the notations
\begin{equation}
\begin{array}{ll}\medskip
\displaystyle{\tilde E_{\alpha\beta}=H_{\alpha\beta}-{H\over 2}h_{\alpha
\beta}},&\tilde E^\alpha_\beta=h^{\alpha\mu}\tilde E_{\mu\beta},\\\medskip
\displaystyle{\tilde E_{ij}=R_{ij}-{R\over 2}g_{ij},}&\tilde  E^i_j=g^{im}
\tilde E_{mj},\\
\displaystyle{\tilde E^{(\alpha)(\beta)}_{(i)(j)}=S^{(\alpha)(\beta)}_{(i)(j)}
-{S\over 2}h^{\alpha\beta}g_{ij}},&\tilde E^{(i)(\beta)}_{(\alpha)(j)}=g^{im}
h_{\alpha\mu}\tilde E^{(\mu)(\beta)}_{(m)(j)}.
\end{array}
\end{equation}

\begin{prop}
The following Einstein d-tensors identities are true:
\begin{equation}
\begin{array}{l}\medskip
\tilde  E^\mu_{\beta/\mu}=0,\\\medskip
\displaystyle{\tilde E^m_{i\vert m}=R^{(m)}_{(\mu)il}P^{l(\mu)}_{\;(m)}-
{1\over 2}g^{kp}R^{(m)}_{(\mu)kl}P^{l\;\;(\mu)}_{pi(m)},}\\
\displaystyle{\tilde E^{(m)(\alpha)}_{(\mu)(i)}\vert^{(\mu)}_{(m)}=S^{(m)(\alpha)
(\delta)}_{(\mu)(i)(l)}S^{(l)(\mu)}_{(\delta)(m)}-{1\over 2}g^{kp}h_{\delta\gamma}
S^{(m)(\gamma)(\delta)}_{(\mu)(k)(l)}S^{l(\alpha)(\mu)}_{p(i)(m)}},
\end{array}
\end{equation}
where $P^{i(\beta)}_{\;(j)}=g^{lm}P^{i\;\;(\beta)}_{lm(j)}$ and
$S^{(i)(\beta)}_{(\alpha)(j)}=g^{lm}h_{\alpha\mu}S^{i(\mu)(\beta)}_{l(m)(j)}$.
\end{prop}
{\bf Proof.} Using the local curvature identities \ref{ci} in the description
of the  d-tensors defined in \ref{dt}, by a direct calculation, we obtain the
following  tensorial identities:
$$
\begin{array}{l}\medskip
H^\alpha_\beta=h^{\mu\varepsilon}H_{\mu\varepsilon\beta}^\alpha,\;\;
R^i_\beta=-g^{lm}R^i_{l\beta m},\;\;R^i_j=g^{ml}R^i_{mlj},\;\;
P^{i(\beta)}_{\;\;(j)}=g^{lm}P_{lm(j)}^{i\;\;(\beta)},\\
P^{(i)}_{(\alpha)\beta}=-g^{lm}h_{\alpha\mu}P^{i\;\;(\mu)}_{l\beta(m)},\;\!
P^{(i)}_{(\alpha)j}=-g^{lm}h_{\alpha\mu}P^{i\;(\mu)}_{lj(m)},\;\!
S^{(i)(\beta)}_{(\alpha)(j)}=g^{lm}h_{\alpha\mu}S^{i(\mu)(\beta)}_{l(m)(j)}.
\end{array}
$$

Let  us consider the following Biachi identities of Cartan canonical connection
\cite{11}
$$
\begin{array}{ll}\medskip
b_1)&\displaystyle{\sum_{\{\alpha,\beta,\gamma\}}H^\delta_{\varepsilon\alpha\beta
/\gamma}=0,}
\\\medskip
b_2)&\displaystyle{\sum_{\{i,j,k\}}\left\{R^l_{pij\vert k}-R^{(m)}_{(\mu)ij}P^{l\;(\mu)}_
{pk(m)}\right\}=0,}\\
b_3)&\displaystyle{
\sum_{\left\{{(\alpha)\atop(i)},{(\beta)\atop(j)},{(\gamma)\atop(k)}\right\}}
\left\{S^{l(\alpha)(\beta)}_{p(i)(j)}\vert^{(\gamma)}_{(k)}+
S^{(m)(\alpha)(\beta)}_{(\mu)(i)(j)}S^{l(\gamma)(\mu)}_{p(k)(m)}\right\}=0.}
\end{array}
$$

Now,  doing the  contractions, $\delta=\beta$ in $b_1$, and $l=j$ in $b_2$,
respectively $b_3$, it follows the Einstein d-tensor  identities required.
q. e. d.

Finally, using the previous proposition and the  relations \ref{*} and \ref{**},
the old conservation laws of Einstein equations, imply
\begin{th}
The new stress-energy d-tensors $\tilde{\cal T}_{\alpha\beta}$, $\tilde{\cal T}_{ij}$
and $\tilde{\cal T}^{(\alpha)(\beta)}_{(i)(j)}$ must verify the following
new conservation laws:
\begin{equation}
\left\{\begin{array}{l}\medskip
\displaystyle{\tilde{\cal T}^\mu_{\beta/\mu}+{1\over 2-n}\tilde{\cal  T}_{M/\beta}
+{1\over 2-pn}\tilde{\cal T}_{v/\beta}=-R^m_{\beta\vert m}-P^{(m)}_{(\mu)\beta}
\vert^{(\mu)}_{(m)}}\\\medskip
\displaystyle{{1\over 2-p}\tilde{\cal T}_{T\vert j}+\tilde{\cal  T}^m_{j\vert m}
+{1\over 2-pn}\tilde{\cal T}_{v\vert j}=-P^{(m)}_{(\mu)j}\vert^{(\mu)}_{(m)}}\\
\displaystyle{{1\over 2-p}\tilde{\cal T}_T\vert^{(\alpha)}_{(i)}+{1\over 2-n}
\tilde{\cal  T}_M\vert^{(\alpha)}_{(i)}+\tilde{\cal T}^{(m)(\alpha)}_{(\mu)(i)}
\vert^{(\mu)}_{(m)}=-P^{m(\alpha)}_{\;\;\;\;i\vert m}},
\end{array}\right.
\end{equation}
\end{th}
where $\tilde{\cal  T}^\alpha_\beta=h^{\alpha\mu}\tilde{\cal  T}_{\mu\beta}$,
$\tilde{\cal  T}^i_j=g^{im}\tilde{\cal  T}_{mj}$ and
$\tilde{\cal  T}^{(i)(\beta)}_{(\alpha)(j)}=g^{im}h_{\alpha\mu}\tilde{\cal
T}^{(\mu)(\beta)}_{(m)(j)}$.
\begin{cor}
If the  local curvature d-tensors $P^{l\;(\mu)}_{pi(m)}$ and $S^{l(\alpha)(\mu)}_
{p(i)(m)}$ vanish, the new conservation laws take the  following simple form:
\begin{equation}
\tilde{\cal T}^\mu_{\beta/\mu}=0,\quad \tilde{\cal T}^m_{i\vert m}=0,\quad
\tilde{\cal T}^{(m)(\alpha)}_{(\mu)(i)}\vert^{(\mu)}_{(m)}=0.
\end{equation}
\end{cor}

\section{Conclusion}

\hspace{5mm} At the end of this paper, we should like to emphasize that the
generalized metrical multi-time Lagrange geometry allows  us to build an
entire theory of physical fields on $J^1(T,M)$, naturally attached to the
following collection of geometrical objects with physical meaning:\medskip

1. A  multi-time Lagrangian function $L:J^1(T,M)\to R$, whose vertical fundamental
metrical d-tensor $G^{(\alpha)(\beta)}_{(i)(j)}(t^\gamma,x^k,x^k_\gamma)$,
regarded as a d-tensor in the indices $(i)$ and $(j)$, is symmetric, of rank
$n$ and having a constant signature;\medskip

2. A fixed semi-Riemannian  metric $h=(h_{\alpha\beta}(t^\gamma))$ on $T$;\medskip

3. A fixed spatial nonlinear connection $N^{(i)}_{(\alpha)j}$ on $J^1(T,M)$.
\medskip

Moreover, if the multi-time Lagrangian function $L$ is a quadratic one (i. e.
the vertical fundamental metrical  d-tensor of $L$ does not depend of partial
directions $x^i_\alpha$), we can construct a natural theory of physical fields,
arised only from $L$ and $h$. This theoretical construction is made, via the
canonical generalized  metrical multi-time Lagrange space $GML^n_p$, attached
to $L$. For more details, see again the example {\bf 2.4}.\medskip\\
{\bf\underline{Open problem}.} The development of an analogous generalized
metrical multi-time Lagrange geometry of physical fields on the jet space of order
two $J^2(T,M)$  is in our attention.

\medskip

\begin{center}
University POLITEHNICA of Bucharest\\
Department of Mathematics I\\
Splaiul Independentei 313\\
77206 Bucharest, Romania\\
e-mail: mircea@mathem.pub.ro\\
\end{center}

\end{document}